\documentclass[12pt]{amsart}

\usepackage{amssymb,amsmath,amsfonts,eurosym,geometry,ulem,graphicx,caption,xcolor,setspace,comment,footmisc,caption,pdflscape,array,hyperref}
\usepackage{enumerate}
\usepackage{subcaption}
\usepackage{graphicx} 

\newtheorem{thm}{Theorem}[section]

\theoremstyle{definition}

\newtheorem{conjecture}[thm]{Conjecture}
\theoremstyle{remark}

\newtheorem{claim}[thm]{Claim}

\newtheorem{question}[thm]{Question}

\DeclareMathOperator{\PSL}{PSL}

\makeatletter
\def\@settitle{\begin{center}%
  \baselineskip14\p@\relax
  \bfseries
  \uppercasenonmath\@title
  \@title
\vskip 0.5cm
     \@subtitle
  \end{center}%
}
\def\subtitle#1{\gdef\@subtitle{#1}}
\def\@subtitle{}
\makeatother

\title{On the topology of the Lorenz system}
\author{Tali Pinsky}
\email{tali@math.tifr.res.in}
\urladdr{http://www.math.tifr.res.in/~tali/}

\begin{document}
\begin{abstract} 
We present a new paradigm for three dimensional chaos, and specifically for the Lorenz equations. The main difficulty in these equations and for a generic flow in dimension three is the existence of singularities. We show how to use knot theory as a way to remove the singularities. Specifically, we claim:\\
(1) For certain parameters, the Lorenz system has an invariant one dimensional curve, which is a trefoil knot. The knot is a union of invariant manifolds of the singular points.\\
(2) The flow is topologically equivalent to an Anosov flow on the complement of this curve, and even to a geodesic flow.  \\
(3) When varying the parameters, the system exhibits topological phase transitions, i.e. for special parameter values, it will be topologically equivalent to an Anosov flow on a knot complement, and different knots appear for different parameter values.

The steps of a mathematical proof of these statements are at different stages. Some have been proven, for some we present numerical evidence and some are still conjectural.
\vspace{0.1in}\\
\vskip -5pt
\noindent\textbf{Keywords:} Lorenz system, knot theory, modular flow.\\
\end{abstract}
\maketitle

\section{Introduction} \label{sec:Introduction}
The Lorenz equations \cite{lorenz1963deterministic}:
\begin{equation}\begin{array}{rl}
\dot{x}(t) & = \sigma(y-x)\\
\dot{y}(t) & = \rho x -y - xz\\
\dot{z}(t) & = xy-\beta z
\end{array}
\end{equation}
originate in weather modeling, but have applications to many other nonlinear phenomena.
They are the principal example of a chaotic system. The parameter values $\rho=28, \sigma=10, \beta=\frac{8}{3}$ were those originally studied by Lorenz and are called the classical values.
At these parameters the Lorenz system possesses the well known butterfly attractor \cite{Tucker1999}.

For $\rho>1$ 
the Lorenz system has three singular points,
one at the origin and two symmetrically related points $p^{\pm}$ at the centers of the butterfly wings.
The singularities prevent the attractor from being hyperbolic 
 \cite{ThreeDimensionalFlows_book,shilnikov2013symbolic}, and are the main reason for the  instability of the system \cite{shilnikov2013symbolic}.
The dynamics slow down near the singularity and this is the main difficulty
in numerical analysis of the flow \cite{araujo2014statistical}.

All three singularities are of saddle type and
have stable and unstable manifolds ~\cite{palis12geom}.
The origin has a one-dimensional unstable manifold, whereas $p^{\pm}$ each have a one-dimensional stable manifold. 
For some parameter values called \emph{T-points}, these manifolds coincide and there are two \emph{heteroclinic orbits}, orbits flowing from $p^\pm$ to the origin.
Numerical studies show that
T-points are central to the dynamics, see for example~\cite{hom10homo,knob14using, shil01meth}, and we show here that this is the case from the topological point of view as well.

We next relate, at each T-point, the Lorenz system to a well known mathematical model for chaos, a \emph{hyperbolic system} \cite{smale}. These are systems with expanding and contracting directions, and although chaotic, can be analyzed and have a well understood statistical behavior. As mentioned, the Lorenz attractor is not hyperbolic,
however we shall see that by removing
the union of the three singular points and their one dimensional invariant manifolds,
the Lorenz flow becomes topologically equivalent to a hyperbolic flow, i.e. there is a continuous invertible map taking orbits to orbits.

\begin{figure}[ht]
  \centering
  \includegraphics[height=4cm]{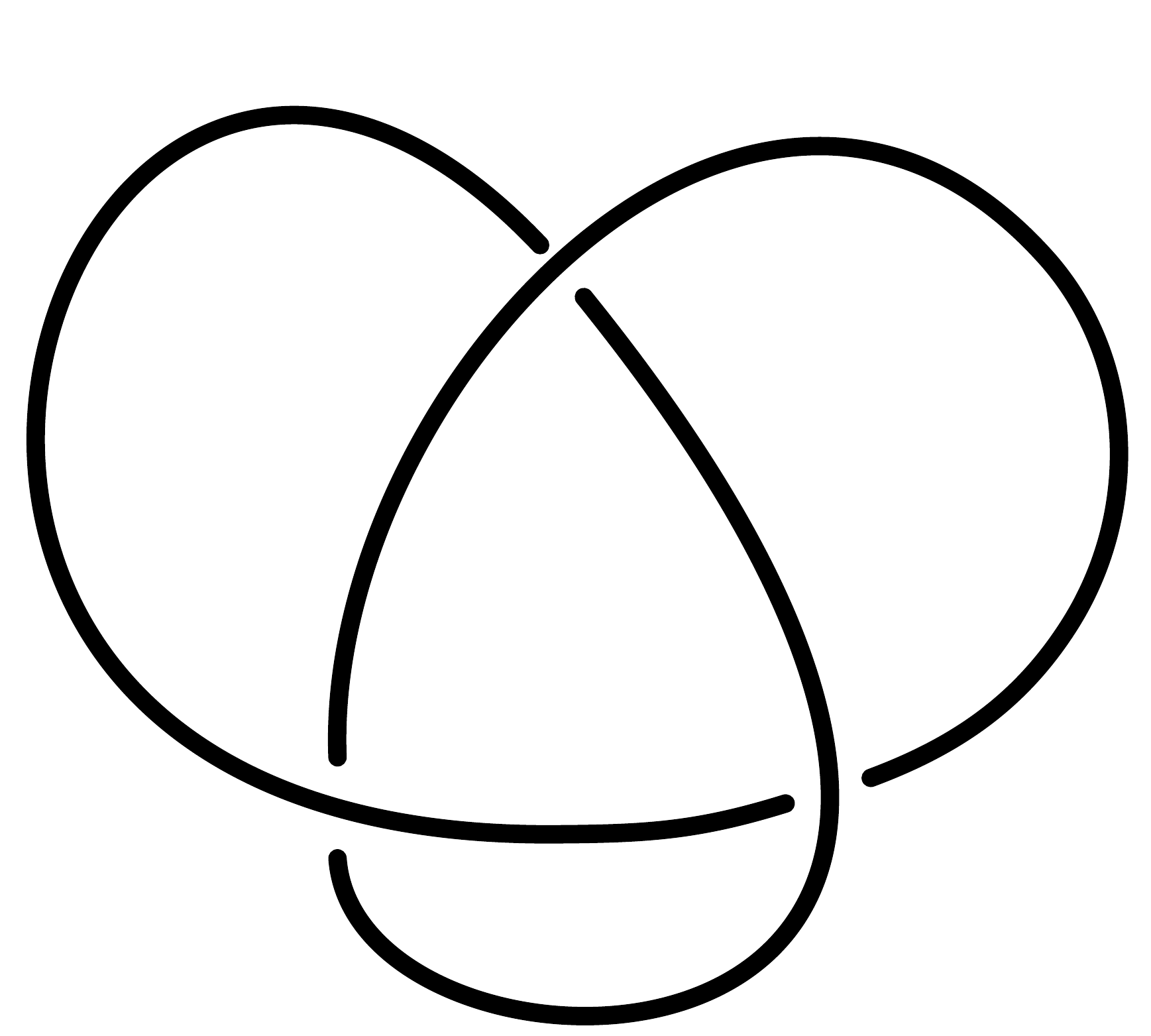}
  \includegraphics[height=4cm]{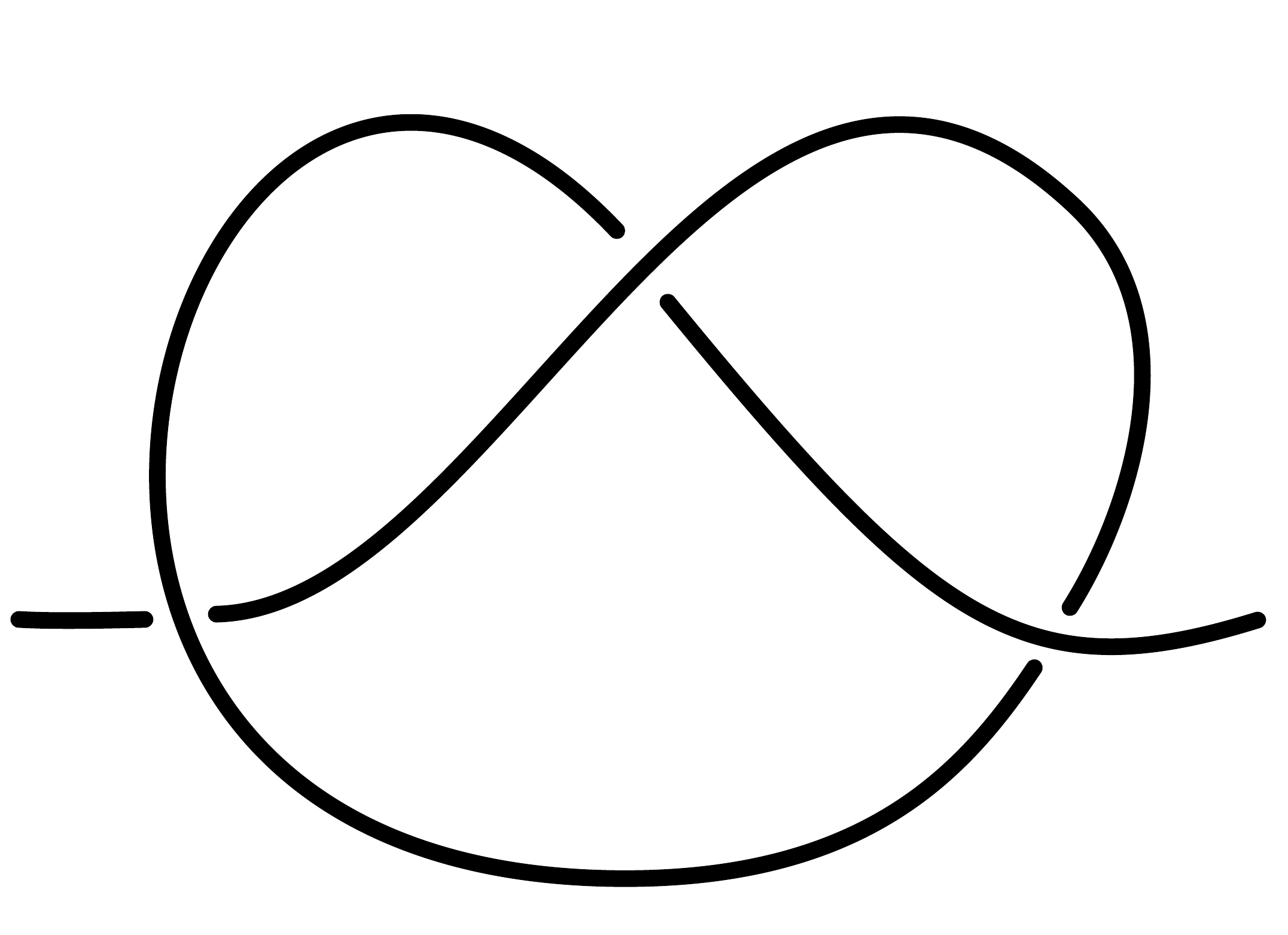}
  \caption{two views of the trefoil knot in $S^3$. The figure shown on the right results from the one on the left when one takes the bottom point of the knot in the figure and pulls it down until it passes through infinity.}
  \label{fig:trefoil}
\end{figure}

A \emph{trefoil knot} is a closed loop
in $S^3$, that can be continuously deformed to the curve given in Figure~\ref{fig:trefoil} without crossing itself.
The first of the parameter values we consider, nearest to the classical parameters, is the primary T-point $\rho_0\approx 30.8680,\,  \sigma_0\approx 10.1673,\, \beta_0=\frac{8}{3}$ ~\cite{alfsen85system}.

\begin{claim}\label{claim:trefoil}
There exists a curve, invariant under Equations (1) for the parameters $\rho_0,\sigma_0,\beta_0$, which is a trefoil knot passing through the three singular points and $\infty$.
\end{claim}

\begin{figure}[ht]
  \centering
    \includegraphics[width=.5\textwidth]{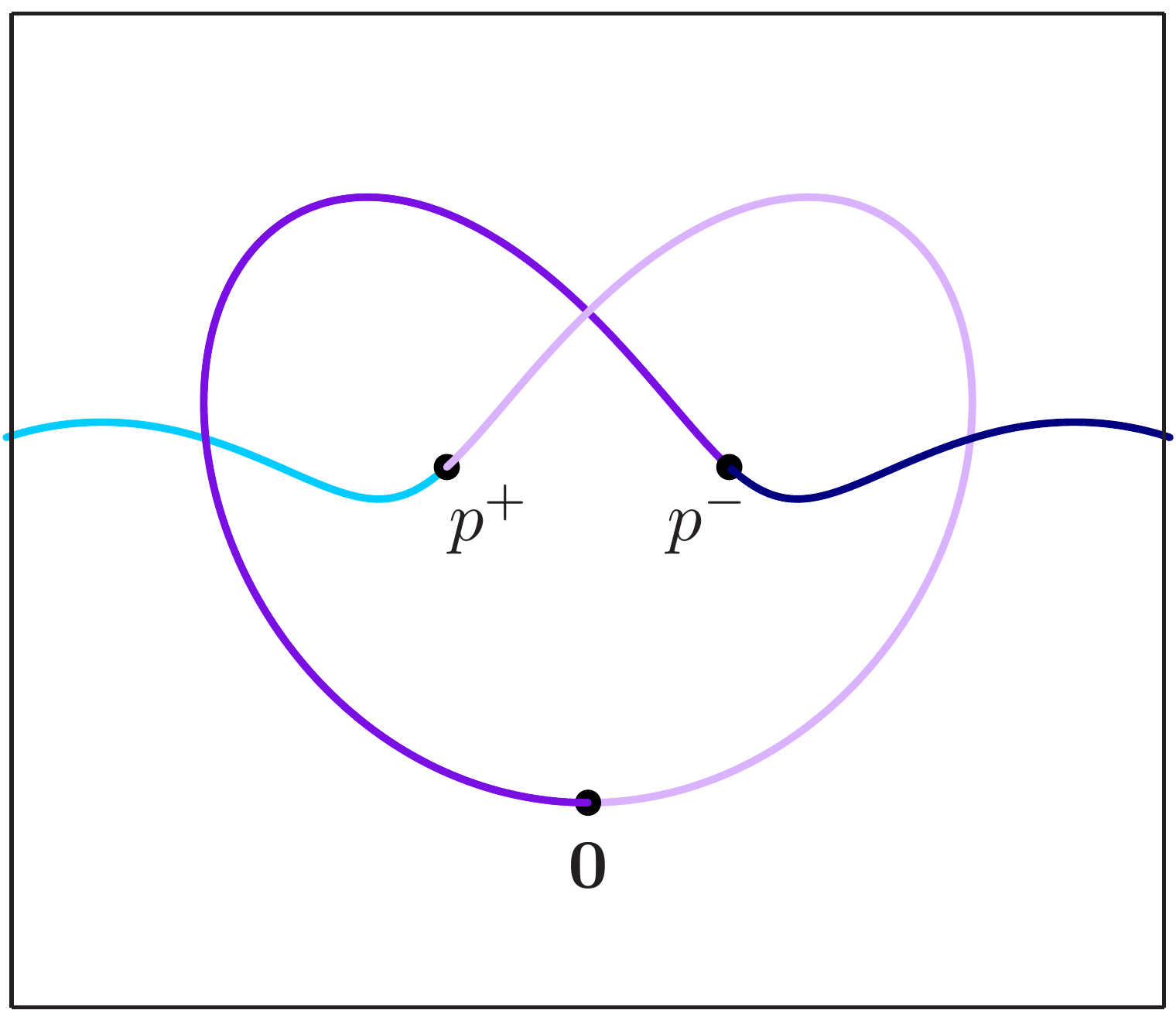}
 \caption{The invariant trefoil for the Lorenz system (Figure courtesy of Jennifer Creaser) \cite{Jen}}
 \label{fig:InvariantTrefoil}
\end{figure}

The invariant trefoil is shown in Figure~\ref{fig:InvariantTrefoil}.
Let us stress that it is not a periodic orbit of the flow. The Lorenz flow has infinitely many periodic orbits with various knot types \cite{BirmanWilliams:KnottedOrbitsI}, and in particular a trefoil shaped periodic orbit. However removing a periodic orbit does not simplify the dynamics of the system, as it does not remove the singularities.

This work was originally motivated by the astonishing similarity between periodic orbits of the Lorenz flow and of the geodesic flow on the modular surface, proven by Ghys \cite{Ghys:KnotsDynamics}. The modular flow is a fundamentally different flow, a mathematical flow originating from number theory. It is defined on the complement of a trefoil knot and thus Claim \ref{claim:trefoil} is key in understanding the relation between the flows.


The next natural question is therefore, could these flows be (in some sense) the same? Even though they have entirely different origins? Even though the Lorenz flow is dissipative, and the modular volume preserving?

We conjecture that, surprisingly, the answer is positive.
In particular, the Lorenz flow is a hyperbolic flow, up to a reparametrization:

\begin{conjecture}\label{conj:reparametrization}
Once removing the invariant trefoil given by Claim ~\ref{claim:trefoil}, the Lorenz flow  is topologically equivalent to the geodesic flow on the modular surface, up to separating the unstable manifold of the cusp.
\end{conjecture}

Essentially, this means that the dissipative nature of the Lorenz equation is entirely due to the fact the point at infinity is repelling. Separating the unstable manifold of the cusp for the modular surface creates the same phenomenon and thus the flows become equivalent.

The trefoil knot shown in figure~\ref{fig:twist} can be considered as the simplest member of a family of knots called \emph{twist knots}. The knot types appearing at the next few T-points are also twist knots and it is natural to expect that all twist knots appear in the Lorenz system for some values of the parameters.
We expect the different knot types and the mechanism of transitioning between them to be key to the global topology of the system.

\section{Numerical method} \label{sec:Numerics}

The trefoil as well as the other knots arising in the system are observed numerically, relying on the computations carried out in \cite{creaser15alpha}, as well as
other studies addressing the existence of the T-points in the Lorenz system, i.e. points where there are two heteroclinic orbits connecting $p^+$ and $p^-$ to the origin.

The first T-point located at $(\rho, \sigma)\approx(30.8680, 10.1673)$ for $\beta=\frac{8}{3}$  was originally found in the 1980's by Petrovskaya and Yudovich~\cite{petrov80homo} and independently by Alfsen and Fr{\o}yland~\cite{alfsen85system}. Thus, the existence of the T-point is well established.


Here, we use the parameter values obtained in \cite{creaser15alpha}, 
The only new ingredient in the present study is that we determine the knot type. To this end we must add two components to the previous results:
\begin{enumerate}[(i)]

\item We compute the other half of the stable manifolds of $p^+$ and $p^-$. 
For a large enough sphere $S$ about the origin, the flow lines crossing it will be directed away from infinity. 
Thus, to show the invariant manifolds connect to infinity we only need to show they connect to such a sphere.

 \item We keep track of the directions of the crossings, that is, which crossing is an undercrossing and which one is an overcrossing. 
Note that the strands in space are actually quite far from each other, and thus a small error in the location of the invariant manifolds cannot change the knot type.
\end{enumerate}

\section{Removing the trefoil knot}
Once we find the one dimensional invariant set as in Claim~\ref{claim:trefoil}, we may remove it from $\mathbb{R}^3$ as its complement is also invariant.
As the curve in question passes through 
the three singular points, defining the flow on its complement results in a non-singular flow which lends itself to classical analysis. 

Once  a knot is removed from $S^3$, the resulting three dimensional space is topologically a \emph{cusped} manifold. That is, the knot itself is ``at infinity''. Moreover, the space can be endowed with a complete metric that reflects its topology.
In this metric, the regular orbits do not slow down near the trefoil, but rather the distances there are large: The distance between two points in this natural metric grows exponentially relative to the Euclidean distance when one approaches the cusp.

Taking into account that the outside of the attracting sphere $S$ 
can be taken as a neighborhood of $\infty$, we choose a tubular neighborhood $N$ of the trefoil knot given in Figure~\ref{fig:InvariantTrefoil}: Let $N$ contain any point that lies outside the sphere $S$, and any point within some distance $\varepsilon$ from the one dimensional manifolds forming the trefoil inside $S$. Topologically $N$ is a solid torus and $\partial N$ is a two dimensional torus.

Let us analyze the behavior of the flow in the neighborhood $N$. We claim $\partial N$ is a union of two transverse annuli:
On an small annulus surrounding the trefoil transversely and centered at the origin, the annulus is close to the stable manifold of the origin. Therefore, the flow enters $N$ near the origin and is transversal to $\partial N$. This corresponds to the fact the unstable manifold of the origin is contained in $N$.
Similarly, in an annulus that bounds the part of $N$ containing both of the points $p^{\pm}$ and the sphere around infinity the flow is escaping $N$.
Choosing $N$ small enough, the two maximal transverse annuli are separated by two meridians of the trefoil.

The precise flowlines within $N$ are determined by the Lyapunov exponents
of the fixed points.
However, its topology is determined by the fact it is in a small neighborhood of the heteroclinic connection.

We claim that once the heteroclinic connection is put at infinity, the flow in $N$ is (almost, in a sense explained in Section \ref{sec:Modular}, topologically equivalent to a flow that is a geodesic flow on a surface of constant negative curvature in a neighborhood of a cusp. Thus, up to changing the metric and the parametrization of the flow, the flow in $N$ is mathematically well understood (see e.g. \cite{Enriquez_StableWindings})
and does not pose the same difficulties as a fixed point. 

\section{Relation to the modular flow}\label{sec:Modular}

\begin{figure}[ht]
  \centering
   \includegraphics[width=.5\textwidth]{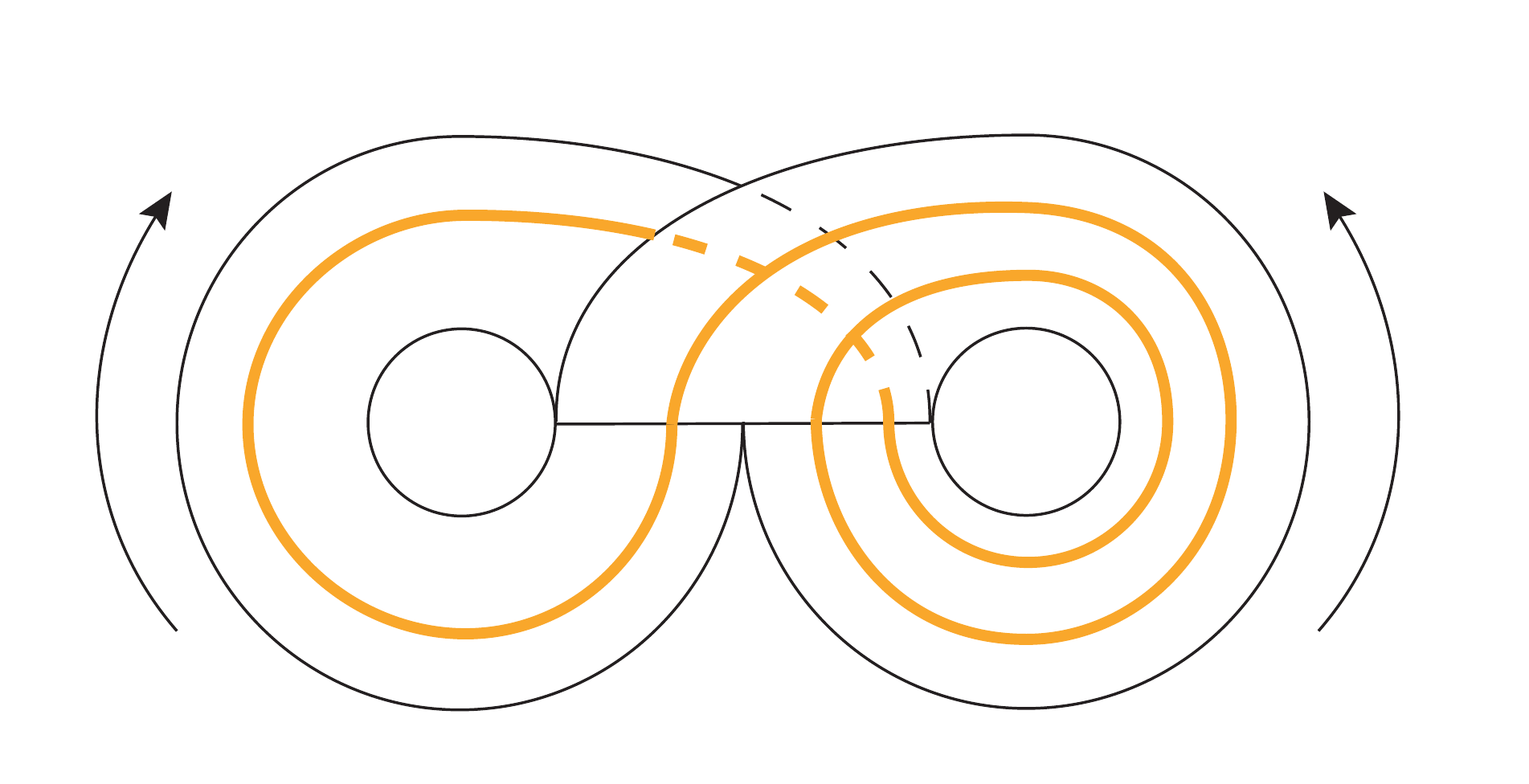}
 \caption{The Lorenz template with a periodic orbit} 
 \label{fig:template}
\end{figure}

For an open set of parameters around the classical ones the recurrent set of the Lorenz equation can be deformed to be contained in a  particular branched surface called the \emph{Lorenz template}, shown in Figure \ref{fig:template}. 
This was proved by Tucker by proving that the Lorenz attractor exists \cite{Tucker1999}. 
The template was used by Birman and Williams  \cite{BirmanWilliams:KnottedOrbitsI} in order to study the topology of the periodic orbits of the flow.

The trefoil complement is topologically equivalent to the matrix group $\PSL_2\mathbb{R}/\PSL_2\mathbb{Z}$. The modular flow is the flow given on this space by left multiplication by the matrix
$$
\left(\begin{array}{cc}
e^t & 0 \\
0 & e^{-t}
\end{array}\right).
$$
This is the geodesic flow on the \emph{modular surface}, the surface obtained by the action of the group $\PSL_2\mathbb{Z}$ on the hyperbolic plane.
A template to this flow exists by \cite{BirmanWilliams2}.

In \cite{Ghys:KnotsDynamics}, Ghys proved that the
template for the modular flow is identical to the Lorenz template.
Thus, the periodic orbits of these two flows are identical. This is surprising as there is no known connection between the two systems. The modular flow is the best possible mathematical flow, not only hyperbolic but manifesting many connections to number theory. See, for example,  \cite{einsiedler2012distribution}, \cite{gurevich2001arithmetic} and \cite{pollicott1986distribution}.

Claim \ref{claim:trefoil} addresses the fact these two flows are not defined on the same space and so is a first step in understanding the relation between the flows.

\begin{figure}[ht]
  \centering
   \includegraphics[width=.5\textwidth]{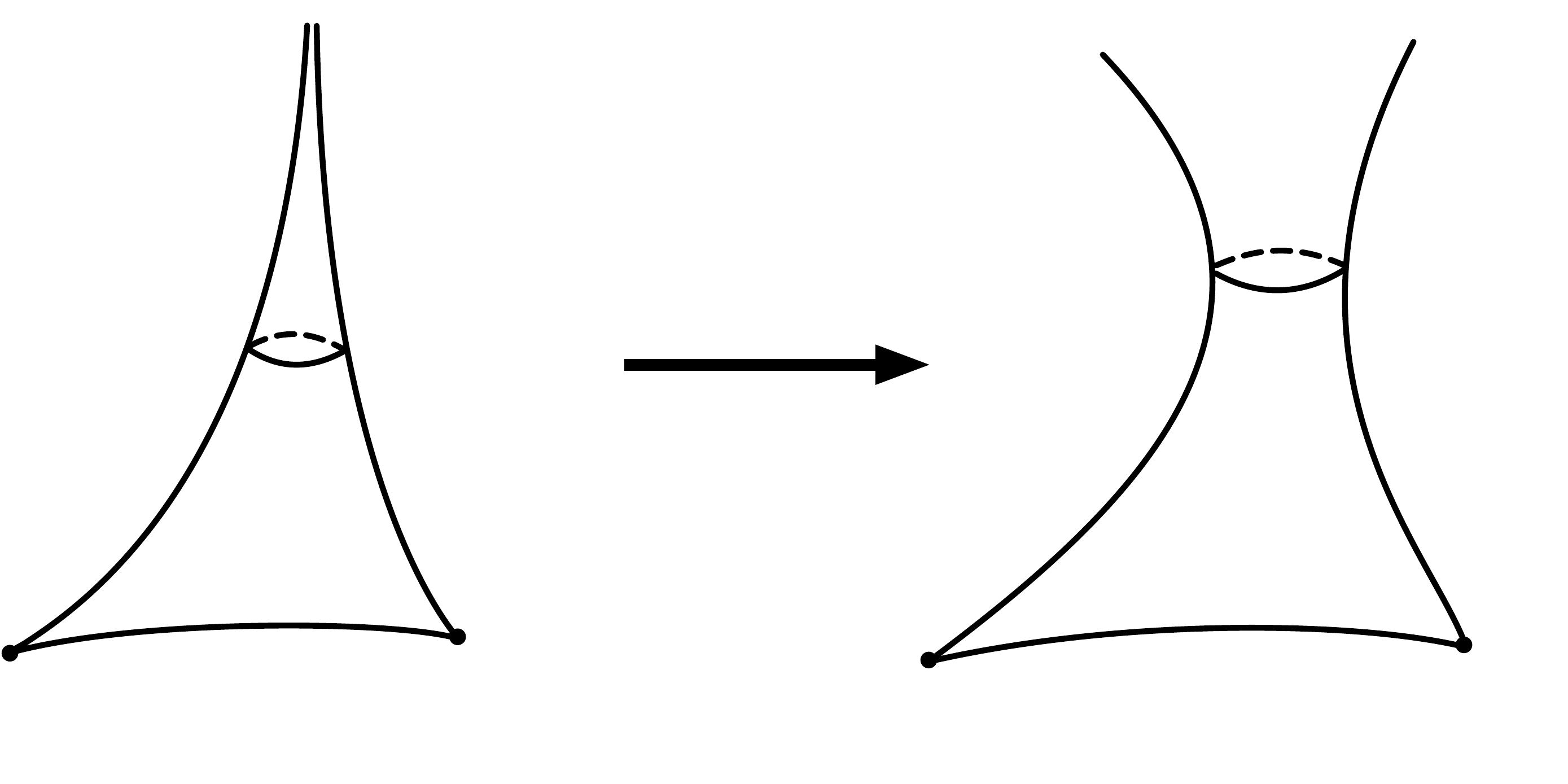}
 \caption{Opening up the cusp for the modular surface.}
 \label{fig:funnel}
\end{figure}

The modular flow is a limit set of a family of geodesic flows on compact hyperbolic surfaces (with boundary). These are obtained by deforming the representation of $\PSL_2(\mathbb{Z})$ so that the cusp becomes an open funnel (see Figure \ref{fig:funnel}), and then truncating them at the unique closed geodesic encircling the cusp.

\begin{figure}[ht]
  \centering
    \includegraphics[width=.7\textwidth]{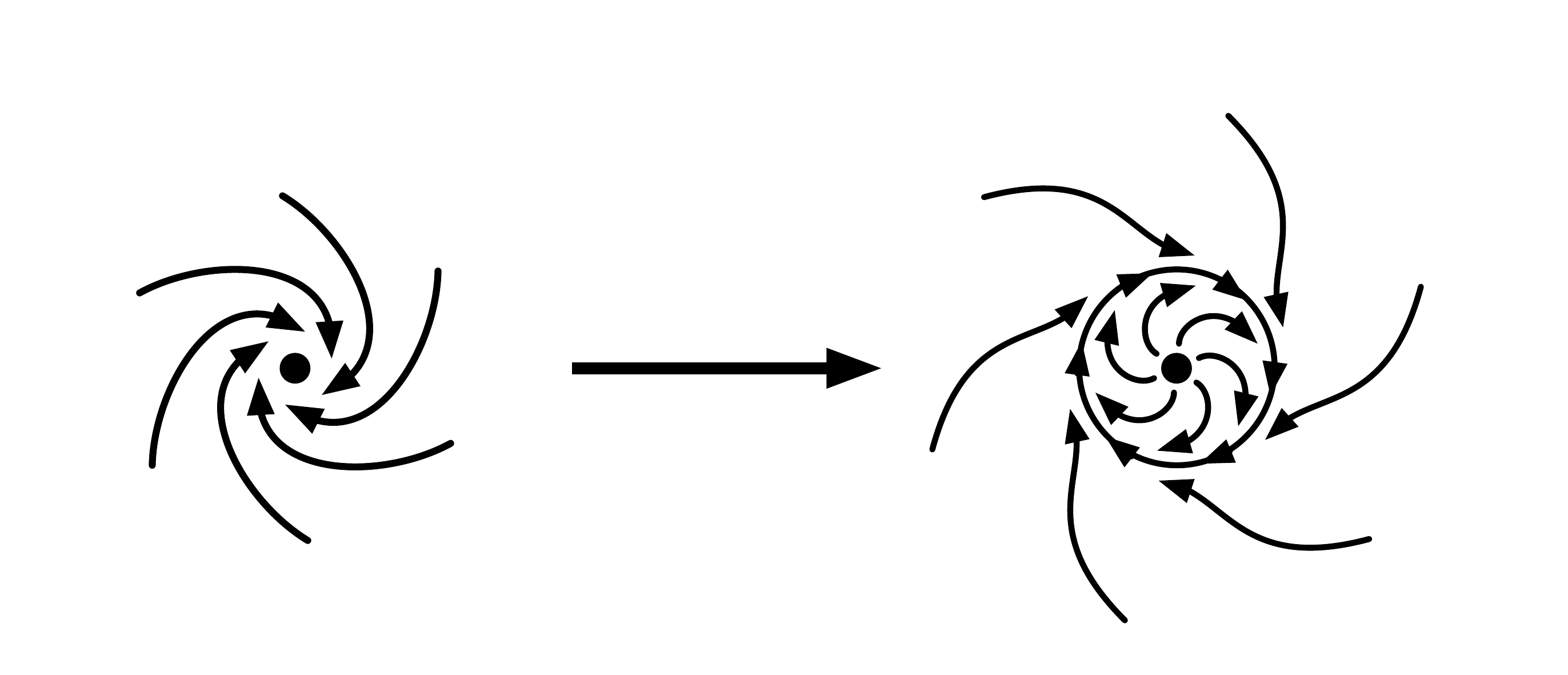}
 \caption{A Hopf bifuraction}
 \label{fig:Hopf}
\end{figure}

The next cornerstone in the proof of Conjecture \ref{conj:reparametrization} is that one may define a blow-up:
Starting from the Lorenz flow at the $T$ point one performes a Hopf bifurcation at each of the wing centers. 
This transforms each singularity to a sink and create two additional orbits, around each singularity.
Such a bifurcation is shown in Figure \ref{fig:Hopf}.

Note this does not affect the trefoil and may be performed in a small neighborhood of any heteroclinic connection of this type.
Next, one may remove a solid torus which is a neighborhood of the trefoil, so that it passes through the two orbits created by the bifurcation, and otherwise is transverse to the flow. Thus the periodic orbits decompose it into two annuli, which the flow crosses in opposite directions. Such a torus is called a Birkhoff torus for the flow.
Note the blow-up has a natural parameter, the diameter of the newly created periodic orbits.

This results in a flow on the compact manifold which is the trefoil knot with toral boundary, just as the usual blow-up \cite{Fried1982CrossSections} (but the boundary torus is not invariant under the flow).

For any small enough blow-up parameter we believe that the blown-up Lorenz flow can be proven to be Axiom A. It will be topologically expanding, as can be seen by its action on a cross section lying in the part of the space between the wing centers. 
We stress that unlike the cross section for the Lorenz flow without the blow up, the return map for the cross section in this case has bounded return times. This corresponds to the fact that the Hopf blow up  separates the recurrent set from the singular points.
This Axiom A flow has a unique (hyperbolic) attractor, which is equivalent to the modular template.


The behavior of the geodesic flow with the opened cusp on the boundary, is identical to the behavior of the Lorenz system on the boundary created by the Hopf blow-up. Namely, the boundary consists of an annulus with the flow flowing inward and an annulus flowing outwards, and these are separated by two \emph{meridional} curves in $S^3\setminus\text{trefoil}$, that are periodic orbits for the flow.
Further, any Anosov flow on the trefoil complement is a the geodesic flow on the modular surface with an open cusp, up to reparametrization.
The identical attractors and boundary behavior are enough to prove:

\begin{claim}
The blown-up Lorenz flow, for any small enough parameter, is orbit equivalent to the geodesic flow on the modular surface with the opened-up cusp.
\end{claim}

\begin{figure}[ht]
  \centering
    \includegraphics[width=.4\textwidth]{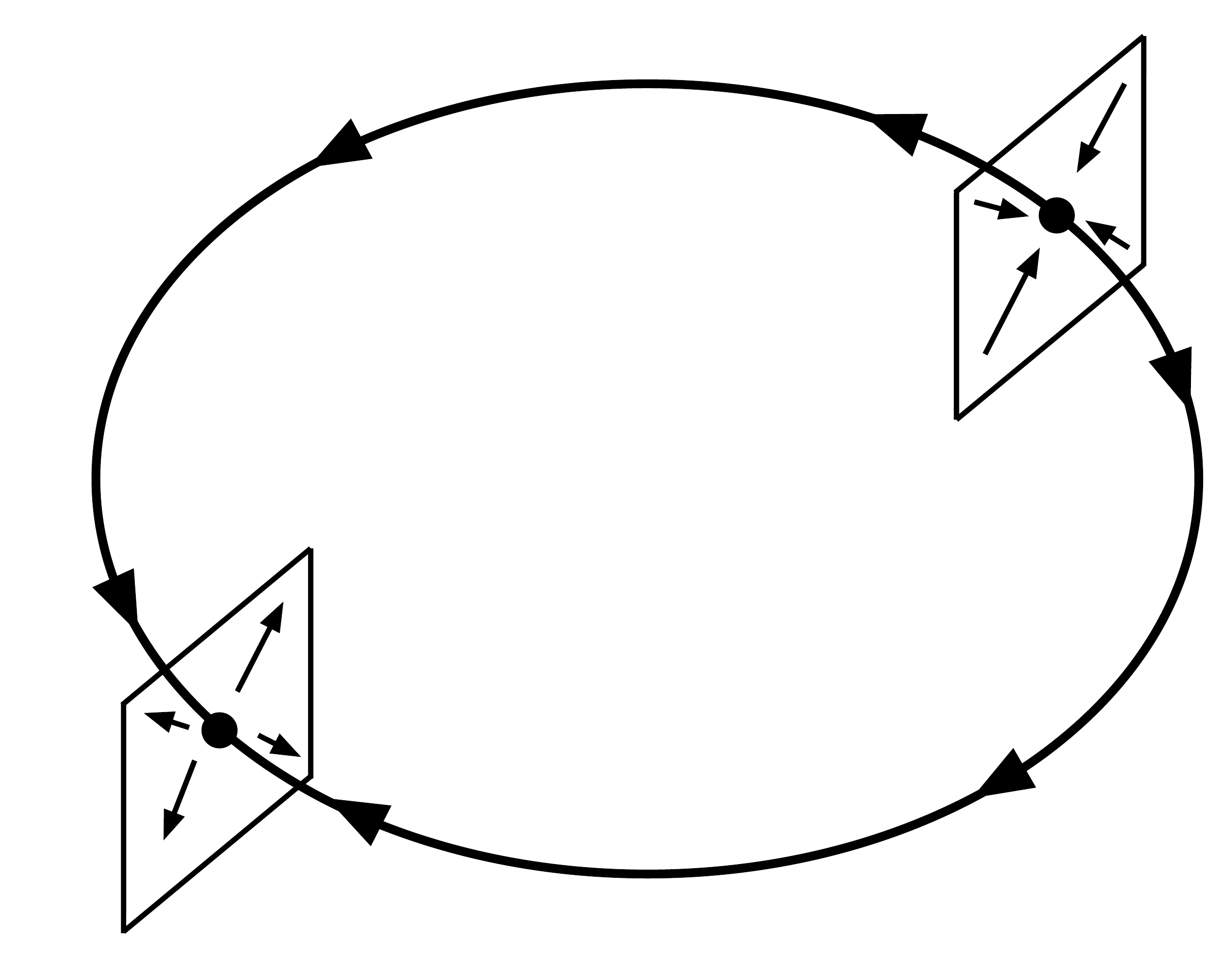}
 \caption{The natural heteroclinic cycle in the completion of the modular flow (when it is embedded as a trefoil knot in $S^3$).}
 \label{fig:Connection}
\end{figure}

Finally, it is natural to expect that one can take the limit when the blow-up parameter goes to zero.
The limit on one side is the Lorenz flow  when removing the one dimensional trefoil without blowing it up, and on the other side it is the modular flow.
Alas, these two flows are too different to be conjugate. One way to see this, is that in the natural completion of the modular flow to a flow on $S^3$, there would be only two singular points, of saddle type.
A neighborhood of the heteroclinic connection connecting these points is Shown in Figure \ref{fig:Connection}. 

\begin{figure}[ht]
  \centering
    \includegraphics[width=.4\textwidth]{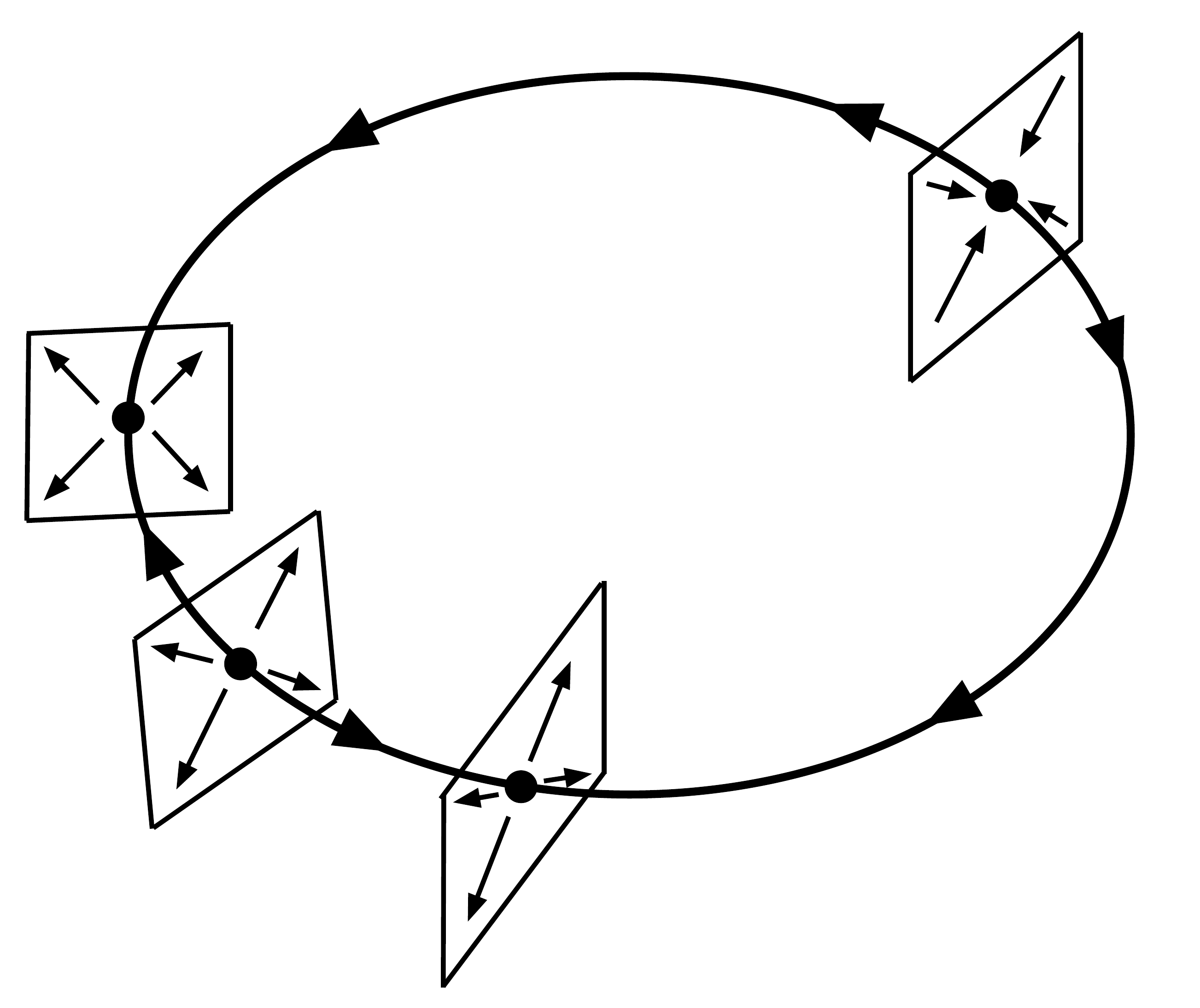}
 \caption{The heteroclinic cycle in the Lorenz flow, consisting of four singular points (including the point at infinity).}
 \label{fig:Connection2}
\end{figure}

On the other hand, The heteroclinic cycle for the Lorenz system is shown in Figure \ref{fig:Connection2}. One of the saddle points is split into two saddles, which are the two wing centers, and a repelling fixed point is added between the two at infinity. 
 
This is what we call separating the cusp's unstable manifold, as when one removes the one dimensional curve from $S^3$ and turns it into a cusp, or a point at infinity, it has a three dimensional part of the space that is its unstable manifold (i.e., that converges to this point as the time goes to $-\infty$), instead of a two dimensional unstable manifold as for the modular flow.

This change turns the volume preserving modular flow into volume decreasing.  In the same way, one could collapse the wing centers onto infinity and transform the Lorenz flow into a volume preserving flow. We claim this is the only difference between the flows and this  would complete the proof of Conjecture \ref{conj:reparametrization}.

We remark that the limit case, the geodesic flow on the cusped modular surface, is indeed unstable. Thus this theory does not contradict the complexity present in the Lorenz system around the T-point.


\section{Varying the parameters \label{sec:bifurcations}}
\begin{figure}[ht]
  \centering
  \begin{subfigure}{5cm}
  \includegraphics[height=4cm]{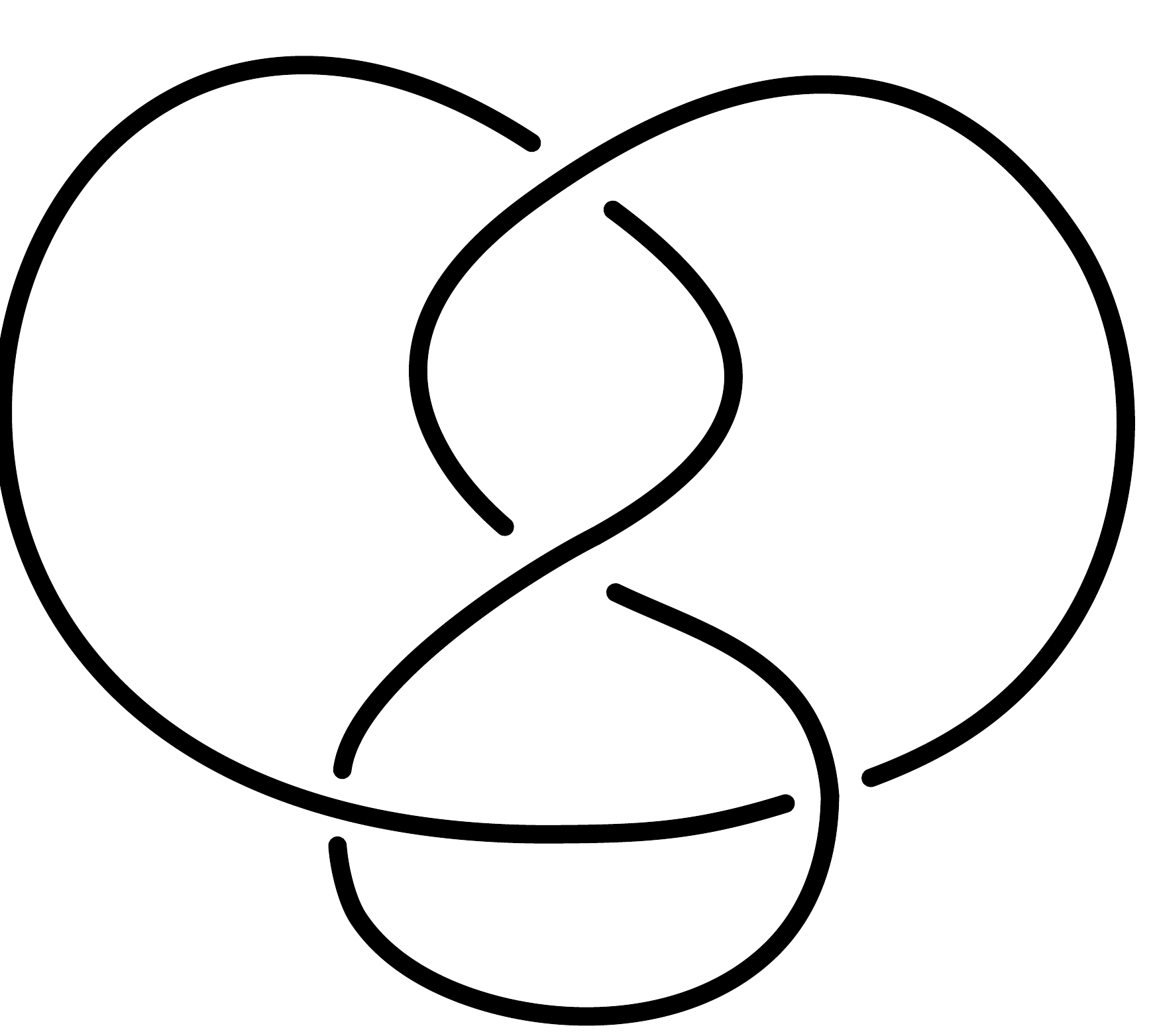}\hspace{1cm}
  \caption{The figure eight knot}
  \end{subfigure}
  \begin{subfigure}{5cm}
  \includegraphics[height=4cm]{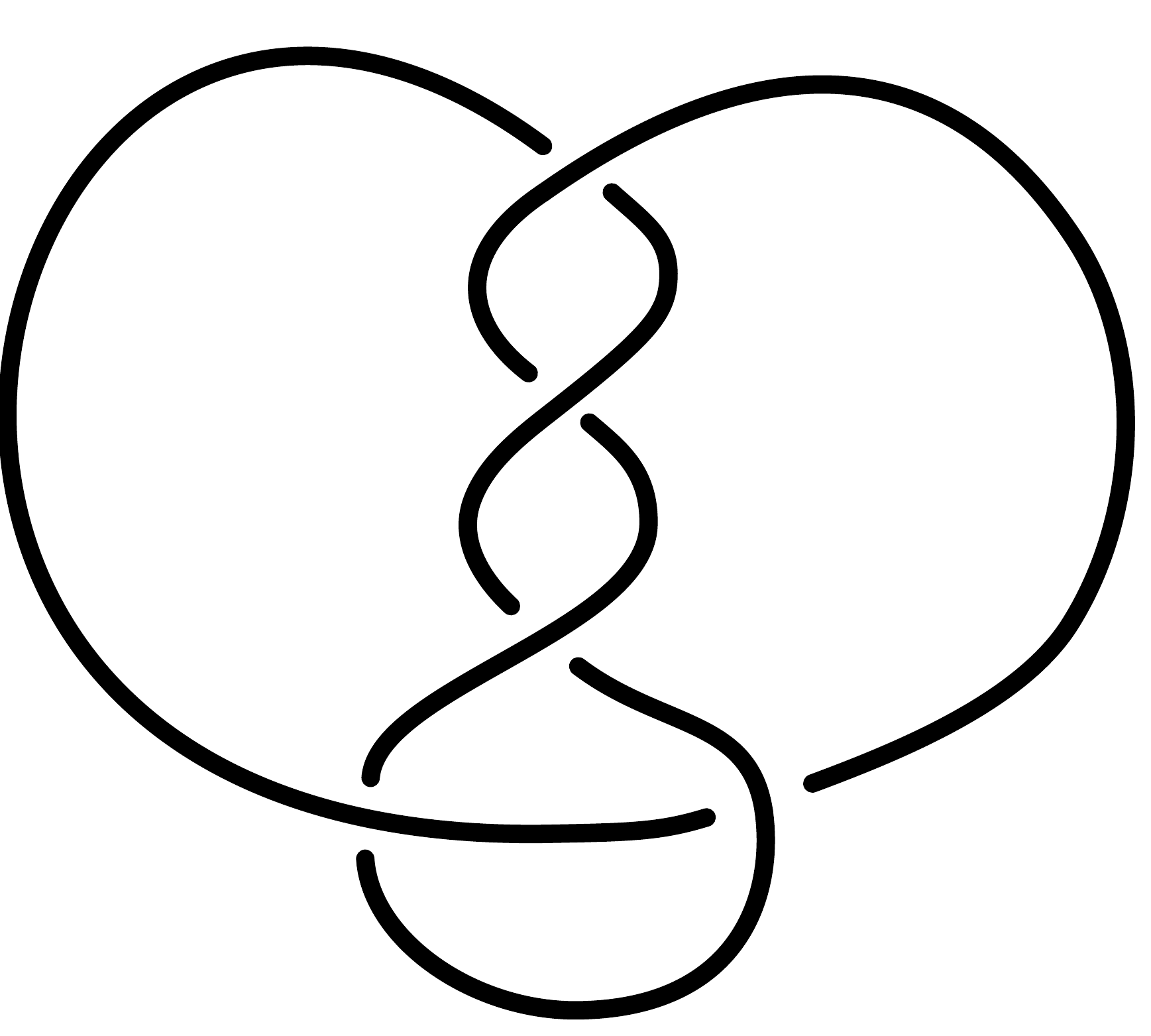}
  \caption{The $5_2$ knot}
  \end{subfigure}
  \caption{the two next members in the family of twist knots. On the left the knot with 2 half twists, called the figure eight knot, and on the right the knot with 3 half twists.}
  \label{fig:twist}
\end{figure}

\begin{figure}[ht]
  \centering
    \begin{subfigure}{6cm}
  \includegraphics[width=5cm]{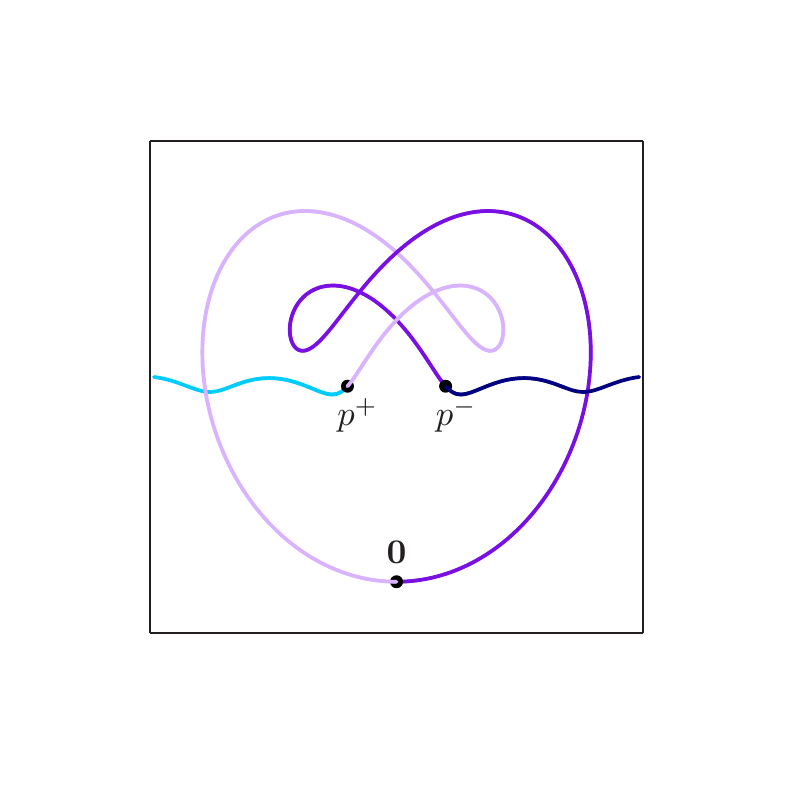}
  \caption{}
    \end{subfigure}
    \begin{subfigure}{6cm}
  \includegraphics[width=5cm]{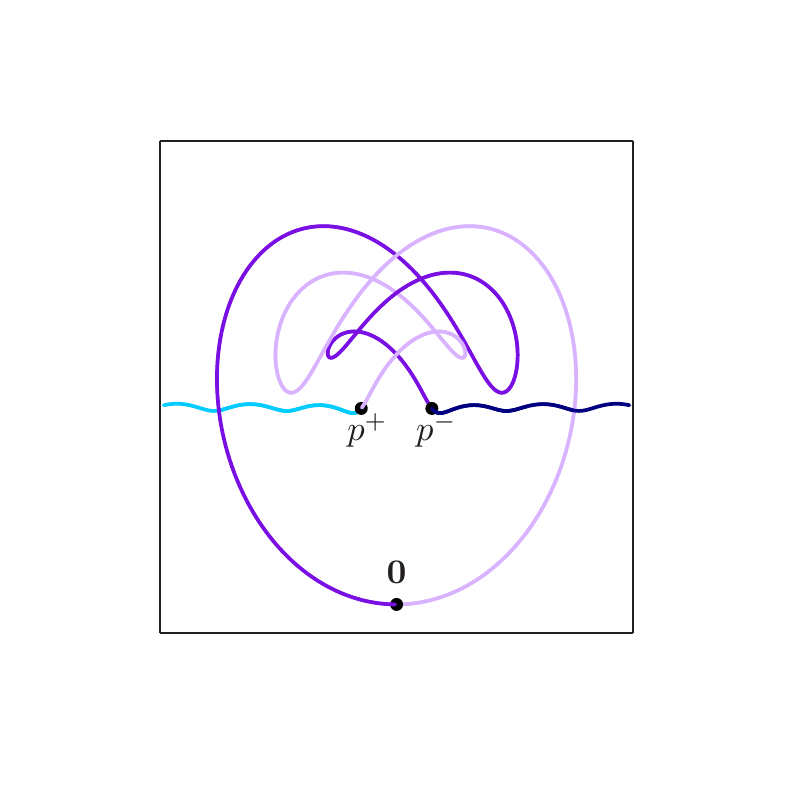}
  \caption{}
    \end{subfigure}
  \caption{The figure eight knot and the $5_2$ knot as invariant knots in the Lorenz system at the next two T-points $(\rho, \sigma)\approx(85.0292,11.8279)$ and $(\rho,\sigma)\approx(164.1376, 12.9661)$ (Figure courtesy of Jennifer Creaser \cite{Jen}).}
  \label{fig:NextKnots}
\end{figure}

As mentioned in the introduction, when one varies the parameters in the Lorenz equations, different knot types appear whenever one has a T-point in the parameter space. 
The next knot one encounters is a knot called the figure eight knot, shown in Figure \ref{fig:twist}(A).
The invariant figure eight knot for the Lorenz system appears for the parameters $(\rho, \sigma)\approx(85.0292,11.8279)$ and is shown in Figure \ref{fig:NextKnots}(A).

The same considerations 
that lead to simpler dynamics in the case of the trefoil knot
will hold for the figure eight knot, as it contains the fixed points and the dynamics on a neighborhood of the knot retains the topology of a geodesic flow about a cusp. 
Hence, the flow on the complement is topologically equivalent to an Anosov flow as well. Interestingly, in contrast to the modular flow, this is a completely new flow that was never studied.

From a knot theoretic viewpoint the trefoil and the figure eight  knot are fundamentally different, the figure eight knot is a  \emph{hyperbolic knot}, that is, its complement has a hyperbolic geometric structure. 
Thus, although the Lorenz flow will be equivalent to an Anosov flow in the complement of the figure eight knot, this flow can never be a geodesic flow and we believe this transition is significant.

The trefoil knot is the first in an infinite family of knots called twist knots.
These are knots obtained by twisting a closed loop, adding to it $n$-half twists, and then cutting it open on one side and re-adjoining the two ends together so that they clasp the loop on the other side. The trefoil knot corresponds to $n=1$, while the knots for $n=2$ (the figure eight knot) and $n=3$ are shown in Figure \ref{fig:twist}.

When $\rho$ and $\sigma$ are increased for fixed $\beta=\frac{8}{3}$ one encounters twist knots with more and more half twists in the Lorenz system. 
The second and third T-point are shown in Figure~\ref{fig:NextKnots}, and it seems the other knots arising for the T-points in \cite[table 3]{creaser15alpha} are all twist knots.

\section{Conclusion \label{sec:conclusions}}

The invariant knots seem to exist along a one dimensional curve \cite{creaser15alpha}. It is indeed expected that the set 
on which the invariant manifolds of $p^\pm$ hit the origin is a set of measure zero.
Nevertheless, these sets are of significance for the dynamics around them and it is intriguing to ask what is the topological explanation for this phenomenon.

\begin{figure}[ht]
\centering
\includegraphics[width=10cm]{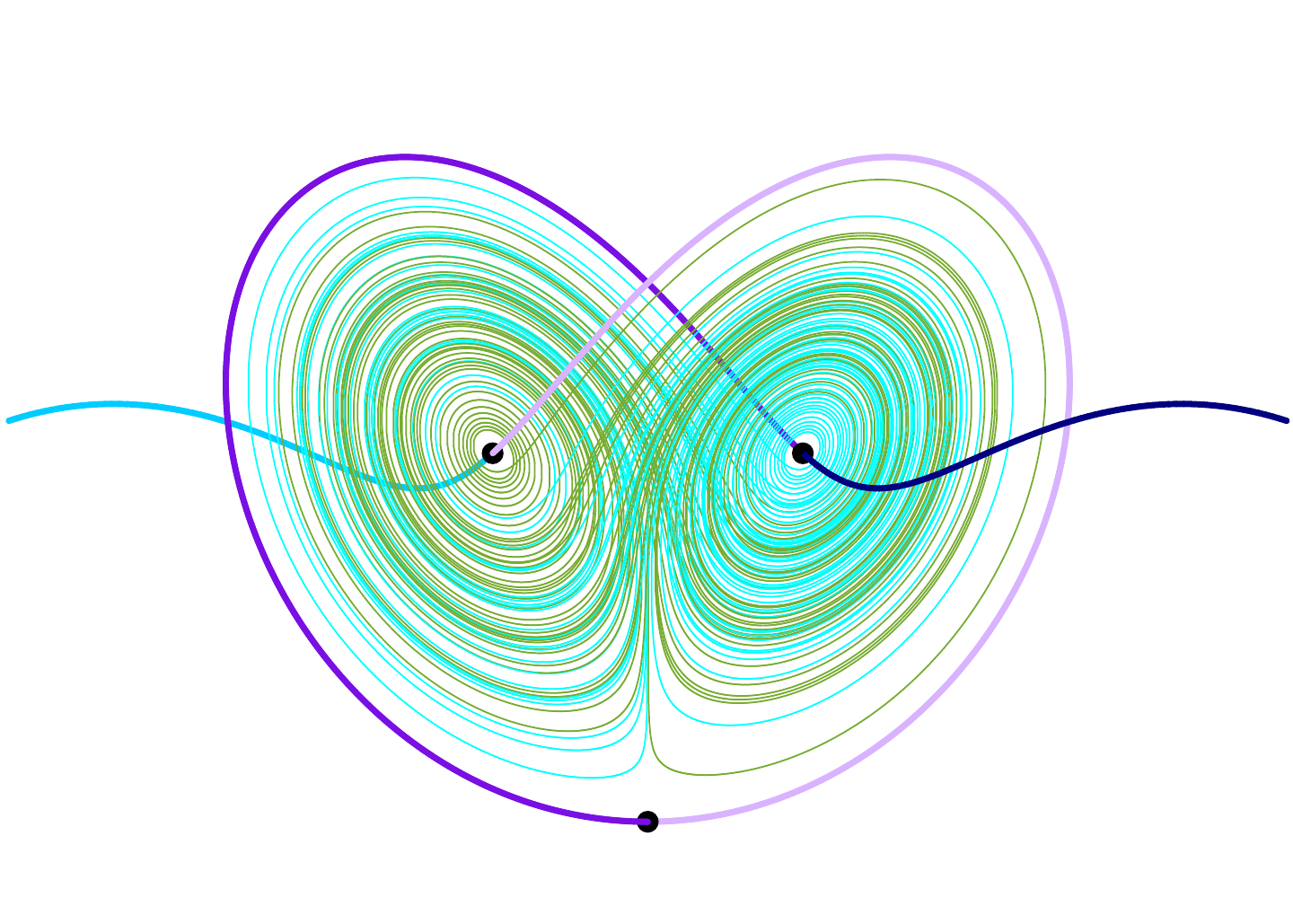}
\caption{The invariant trefoil is nearly determined by the fact it is the boundary of the Lorenz butterfly.  }\label{fig:trefoil_attractor} 
\end{figure}

In Figure ~\ref{fig:trefoil_attractor} we depict the same invariant trefoil as in Figure ~\ref{fig:InvariantTrefoil}, together with the attractor. It can be seen that part of the trefoil (the heteroclinic connections) is the boundary of the attractor 
(c.f. \cite[pages 36-37]{Sparrow}). 

\begin{figure}[ht]
\centering
\includegraphics[width=7cm]{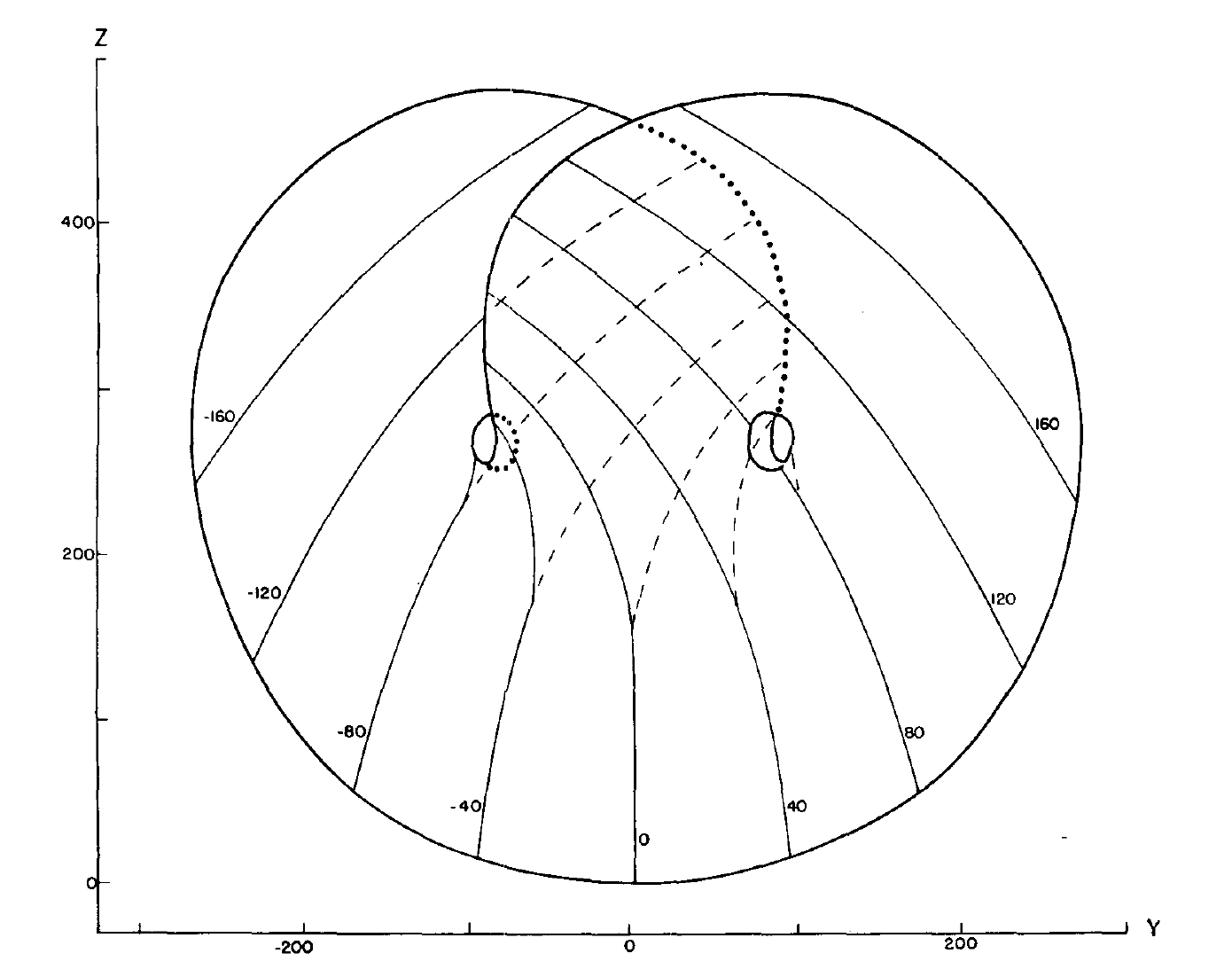}
\caption{Lorenz' diagram of the attractor. The boundary lies on a thick trefoil knot.}\label{fig:LorenzDiagram} 
\end{figure}

On  some parameters off the curve of existence of the trefoil, the attractor still exists, with the same topology and (topologically) the same boundary curve. In Lorenz' original paper, he discusses the fact the attractor can be approximated by two 2-dimensional bands, coming together along the segment connecting $p^+$ to $p^-$. Lorenz' diagram of this approximated attractor is shown in Figure \ref{fig:LorenzDiagram}. The boundary of this two dimensional model of the attractor is part of the unstable manifold of the origin. It is very similar to the two heteroclinic connections that are part of the trefoil. However, for the classical parameters there is no heteroclinic connection, the unstable manifold of the origin does not connect to the other fixed points. 

It is compelling to ask whether one could remove a thickened version of the invariant trefoil from the space, on which the boundary of the attractor lies (c.f. \cite[Figure 4]{barrio12knead}):
Although the unstable manifold of the origin misses the other two fixed points, it misses them only slightly. Thus, a small ball around the origin would hit both wing centers, and then continue to infinity on their other side.

It seems that the flow then becomes another well understood mathematical model, the geodesic flow on a truncated modular surface (see \cite{GhysJapanese}). Thus for an open set of parameters in the Lorenz equation, the flow is equivalent to the geodesic flow on the modular surface with different parameters.

This gives a new hierarchy of the periodic orbit sets (and the entire recurrent dynamics). For this set of parameters, The curve on which the trefoil connection exists is the curve on which the system has the largest set of periodic orbits. 
Furthermore, the further from this curve the system is, the more truncated it is, and the smaller its set of periodic orbits.

Another important question from the topological viewpoint is the following.

\begin{question}
What is the topological mechanism of transitioning between the trefoil and the figure eight knot?
\end{question}

There seem to be two possible explanations:
\begin{enumerate}
\item When the stable manifolds of $p^{\pm}$ miss the origin they continue to infinity on their other side as well. Thus the complement of the three invariant manifolds becomes the complement of a bouquet of three circles passing through the point at infinity. 
The flow on their complement has two different limit points where the space reduces to a knot complement. One is the modular flow on the trefoil complement and the other a flow on the figure eight knot complement.

\item There is an invariant \emph{tunnel}, an  arc with both its endpoints on the knot so that  these two knots together with the tunnel become equivalent. This would imply that in fact the three dimensional manifold is fixed, and the Lorenz flow is topologically equivalent to the modular flow throughout all these different points in the parameter space. 
\end{enumerate}

From the point of view of dynamical systems, an obvious question is in what way can the topological equivalence be used to establish dynamical properties of the system. 

\begin{question}
Can one establish an exponential decay of correlations or a central limit theorem on the trefoil complement, in a cusped metric?
\end{question}

This is strongly related to the (completely open) question of how well behaved is the orbit equivalence between these two flows.

\subsection*{Acknowledgments}
The author wishes to thank Joan Birman, Christian Bonatti, \'Etienne Ghys and Amos Nevo for their encouragement and for helpful discussions, and is grateful to Jennifer Creaser for generously sharing her knowledge of the Lorenz system, carrying out the numerical simulations mentioned in this paper and generating the figures. This research was supported by a UGC grant.
\bibliographystyle{plain}
\bibliography{bibliography.bib}

\end{document}